\documentclass[11pt]{article}
\usepackage{amsmath,amssymb}
\usepackage{booktabs}
\usepackage{threeparttable} 
\usepackage{geometry}
\usepackage{booktabs}
\usepackage[numbers]{natbib}
\usepackage{authblk}
\usepackage[colorlinks=true, linkcolor=blue, citecolor=blue, urlcolor=blue]{hyperref}
\usepackage{pgfplots}
\usepgfplotslibrary{groupplots}
\pgfplotsset{compat=1.18} 
\title{On the Limits of Interpretable Machine Learning in Quintic Root Classification}

\author{Rohan Thomas, Majid Bani-Yaghoub\footnote{Corresponding author. Email: baniyaghoubm@umkc.edu}}
\affil{\small  
 Division of Computing, Analytics \& Mathematics, School of Science and Engineering, University of Missouri-Kansas City, 5100 Rockhill Rd., Kansas City, Missouri 64110, USA
}
\date{}
\begin{document}
\maketitle

\begin{abstract}

Can Machine Learning (ML) autonomously recover interpretable mathematical structure from raw numerical data? We aim to answer this question using the classification of real-root configurations of polynomials up to degree five as a structured benchmark. While algebraic discriminants exist for degrees two through four, quintic polynomials admit no general symbolic solution in radicals. This provides a controlled setting to evaluate whether ML models can rediscover meaningful mathematical invariants. We tested an extensive set of ML models, including decision trees, logistic regression, support vector machines, random forest, gradient boosting, XGBoost, symbolic regression, and neural networks. Neural networks achieved strong in-distribution performance on quintic classification using raw coefficients alone (84.3\% ± 0.9\% balanced accuracy), whereas decision trees perform substantially worse (59.9\% ± 0.9\%). However, when provided with an explicit feature capturing sign changes at critical points, decision trees match neural performance (84.2\% ± 1.2\%) and yield explicit classification rules. Knowledge distillation reveals that this single invariant accounts for 97.5\% of the extracted decision structure. Out-of-distribution, data-efficiency, and noise robustness analyses indicate that neural networks learn continuous, data-dependent geometric approximations of the decision boundary rather than recovering scale-invariant symbolic rules. This distinction between geometric approximation and symbolic invariance explains the gap between predictive performance and interpretability observed across models. Although high predictive accuracy is attainable, we find no evidence that the evaluated ML models autonomously recover discrete, human-interpretable mathematical rules from raw coefficients. These results suggest that, in structured mathematical domains, interpretability may require explicit structural inductive bias rather than purely data-driven approximation.
\end{abstract}



\section{Introduction}
While Machine Learning  (ML) models have achieved significant predictive accuracy across various domains, enhancing their interpretability remains a persistent challenge. In the present work, we address the following question. Can ML discover interpretable mathematical rules from raw data? We define “autonomous discovery” as the ability of an ML model to (i) learn a decision boundary directly from raw coefficients without mathematically engineered invariants, and (ii) produce an explicit, human-interpretable rule that generalizes beyond the training distribution without additional symbolic guidance. Recent work has showcased ML abilities in mathematical pattern recognition \cite{davies2021,fawzi2022}, but the learned patterns remain entrapped in network weights. Our research investigates whether the rules learned by ML models can help mathematicians uncover discrete, symbolic decision rules governing real-root structure, rather than merely approximate classification boundaries. Quintic polynomial classification offers the ideal testbed. Although the Abel-Ruffini theorem proves no general solution in radicals exists \cite{abel1824,ruffini1799, Stewart2015Galois}, real-root counting remains algorithmically solvable through classical analytic tools (e.g., Sturm sequences). Thus, quintic classification presents a structured but nontrivial benchmark for evaluating whether ML can rediscover meaningful mathematical structure from raw coefficients. Real and complex roots of quintic polynomials can be numerically estimated with high precision, which provides ground truth for the present study. Table \ref{tab:ml_models} lists popular ML models and their descriptions. In the present work, we evaluated a diverse set of machine learning approaches, 
including decision trees \cite{Breiman1984CART, soysal2022machine}, logistic regression, support vector machines 
\cite{CortesVapnik1995SVM}, random forests \cite{Breiman2001RandomForest}, 
gradient boosting machines \cite{Friedman2001GBM}, XGBoost 
\cite{arjmand2025comparative, ChenGuestrin2016XGBoost}, and neural networks 
\cite{Rosenblatt1958Perceptron,Bishop2006PRML}, to investigate this question. 
Despite utilizing both classical statistical learning frameworks and modern ensemble and deep learning architectures, the models did not autonomously discover interpretable mathematical rules as we had initially hypothesized.\\

The main contributions of the present work are as follows. First, we test the interpretability of ML models on degrees 2-4 (where algebraic discriminants exist) and confirm that, as long as they are provided with the necessary features, interpretable models could successfully extract near-perfect classification rules. Second, we demonstrate that neural networks perform well in quintic root classification with only coefficient features (84.2\%), suggesting some form of pattern learning, although the nature of this learning remains unclear. Third, we show that decision trees, when evaluated through the same framework, achieve only 60.7\% balanced accuracy, demonstrating a significant capability gap. Finally, through decision-tree distillation, we identify Crit8 (critical point sign changes) as 97.5\% responsible for learned structure.\\

While the results of our paper yield some promise, the ML models ultimately failed to provide extractable mathematical rules without human intervention. This suggests that autonomous discovery of interpretable mathematical rules remains an open challenge for the scientific community. 
Recent work has investigated the usefulness of AI and ML in mathematics \cite{golladay2025evaluating}. Davies et al. \cite{davies2021} used supervised learning to guide conjecture formation in knot theory, where human mathematicians provided pattern labels. Fawzi et al. \cite{fawzi2022} discovered faster matrix multiplication algorithms through reinforcement learning. Romera-Paredes et al. \cite{romera2024} employed large language models (LLMs) with program search for combinatorial construction. Lample and Charton \cite{lample2020} trained transformers on symbolic mathematics tasks, including integration and differential equations. Although these approaches show promise, they generally work with explicit problem formulations where the true extent of ML innovation is fairly limited. Our work differs by investigating whether models can discover mathematical structure from raw numerical data (polynomial coefficients) without symbolic guidance and, critically, whether the learned knowledge can be extracted.

\section{Methods}

\subsection{Quintic Polynomial Experimental Setup}

For our experiments, we generated 40,000 quintic polynomials with coefficients sampled from $[-10, 10]$. Roots were computed using NumPy's eigenvalue method with the threshold $10^{-10}$ for real/complex classification. Our experiments returned three classes: 5 real roots (class 0), 3 real roots (class 1), and 1 real root (class 2). Model performance was evaluated using 5-Fold Stratified Cross-Validation across 20 independent seeds (0-19). All metrics are reported as the mean $\pm$ a 95\% Confidence Interval. Balanced accuracy was used to mitigate potential bias due to class imbalance.
As explained below, we first evaluated our scaffold over lower-degree polynomials where algebraic discriminants exist to validate that decision trees can extract known mathematical rules when provided with the necessary features.\\

For quadratic polynomials $ax^2 + bx + c = 0$, the discriminant $\Delta = b^2 - 4ac$ algebraically determines root type. $\Delta \geq 0$ indicates real roots while $\Delta < 0$ indicates complex ones. To replicate the decision boundary derived from the quadratic discriminant, we provided neural networks and decision trees with the discriminant ratio feature $b^2/ac$. Additionally, we tested Symbolic Regression (PySR) on quadratic coefficients to determine whether it could bypass the need for the discriminant ratio feature by deriving the true discriminant formula. \\
 
 For monic cubic polynomials $x^3 + Ax^2 + Bx + C$, we utilized reduced form parameters: $\alpha = (A/3)^2 - B/3$ and $\beta = 2(A/3)^3 - (A/3)(B/3) + C$ to determine root type. We then aggregated these parameters into the ratio $\beta^2/\alpha^3$, which serves as the key discriminant feature for distinguishing between the cases of all real roots and one real with two complex conjugates.\\

For quartic polynomials, classical invariants exist but are more complex. So, instead of transcribing them directly, we provided Tschirnhaus-style invariants and discriminant-related features to replicate their decision boundary. For monic quartics $x^4 + Ax^3 + Bx^2 + Cx + D$, relevant invariants include:
\begin{align*}
I &= 12e - 3bd + c^2 \\
J &= 72ace + 9bcd - 27ad^2 - 27b^2e - 2c^3
\end{align*}
where $a=1$, $b=A$, $c=B$, $d=C$, $e=D$.\\


Next, we evaluated an extensive set of ML approaches to identify which model performed best on quintic classification with coefficient features alone. This initial screening, evaluated across three different seeds, would help us determine the model we should optimize with explicit feature engineering

\begin{table}[htbp]
\caption{Overview of Machine Learning models used in this study to evaluate their interpretability}
\centering
\renewcommand{\arraystretch}{1.3}
\begin{tabular}{p{5cm} p{11cm}}
\hline
\textbf{Model} & \textbf{Description} \\
\hline
\textbf{Decision Trees (CART)} & Interpretable models that produce explicit if-then rules, allowing direct examination of learned structure. \\
\textbf{Logistic Regression} & Linear classifier to test whether linear decision boundaries suffice for this problem. \\
\textbf{Support Vector Machines (SVM)} & RBF kernel to explore whether kernel methods could find separating boundaries. \\
\textbf{Random Forest} & Ensemble of decision trees to assess whether aggregation improves performance. \\
\textbf{Gradient Boosting} & Sequential ensemble method for comparison. \\
\textbf{XGBoost} & Optimized gradient boosting implementation. \\
\textbf{Symbolic Regression} & Genetic programming approach (PySR) to search for exact mathematical formulas mapping coefficients to root counts. \\
\textbf{Neural Networks} & Multi-layer perceptrons to explore whether deeper models could learn complex patterns. \\
\hline
\end{tabular}
\label{tab:ml_models}
\end{table}

Neural Networks consistently outperformed other models on raw coefficient inputs. Decision trees, on the other hand, performed substantially worse, revealing the first signs of the complexity gap between powerful models and interpretable ones. Despite their underwhelming performance, decision trees are a necessary foil to the "black box" complexity of Neural Networks, which were evaluated alongside decision trees in a more robust 20-trial evaluation with explicit feature engineering.

\subsection{Mathematical Features for Quintic Classification}

According to the Abel-Ruffini theorem \cite{ Stewart2015Galois}, unlike their second, third, and fourth degree counterparts, quintic polynomials do not have a clean algebraic decision boundary. So, to replicate quintic behavior as closely as possible, we engineered 63 feature coefficients and six classical methods from polynomial theory. In the following, we explain how these classical methods have been utilized in our study. \\


Following Sturm Sequences (8 features), we construct the complete Sturm sequence through iterated polynomial division, starting from $p(x)$ and its derivative \cite{sturm1829}. This sequence allows us to count roots, since the difference in sign changes between two points equals the number of real roots in that interval. Our features extract sign changes at $\pm\infty$ (indicating total real roots), plus sign changes evaluated at five other points $\{-10, -1, 0, 1, 10\}$ to capture local roots.\\


Using Descartes' Rule (6 features), we count sign changes in the coefficient sequences of $p(x)$ and $p(-x)$ to bound positive and negative real roots \cite{descartes1637}. Features include the number of positive root sign changes, negative root sign changes, their sum (total real root bound), the implied minimum complex root count ($5 - \text{total bound}$), and parity indicators for both counts.\\


Utilizing Newton's Sums (10 features), we compute power sums $s_k = \sum_{i=1}^{5} r_i^k$ directly from coefficients using Newton's identities without finding the actual roots \cite{newton1707}. The first five features are the raw power sums $s_1$ through $s_5$. We then derive the estimated mean of roots ($s_1/5$), the estimated variance ($s_2/5 - (s_1/5)^2$), and three normalized ratios that measure properties like how spread out or clustered the roots are in scale-invariant form.\\


 Using Critical Points (10 features), we find the roots of the derivative $p'(x) = 5x^4 + 4Ax^3 + 3Bx^2 + 2Cx + D$ to identify critical points where the polynomial changes from increasing to decreasing (or vice versa) \cite{marden1949}. Features include: the count of real critical points, their position statistics (minimum, maximum, mean, standard deviation), polynomial values at these points (minimum, maximum, mean), and the number of real inflection points. Most importantly, we compute Crit8: the number of sign changes in the sequence $\{p(c_1), p(c_2), \ldots, p(c_k)\}$ where $c_i$ are ordered critical points. By the Intermediate Value Theorem, each sign change between consecutive extrema guarantees a root in that interval, making Crit8 a vital feature in quintic classification. \\


Following Hybrid Symbolic approach (16 features), we compute Tschirnhaus invariants $I_2$ through $I_5$ (generalizations of classical discriminants), four novel algebraic combinations such as $A \cdot B \cdot C - D \cdot E$ and $A^2E - B^2D + C^3$, three discriminant-like difference measures ($B^2 - AC$, $C^2 - BD$, $D^2 - CE$), and five scale-invariant coefficient ratios. These features capture algebraic relationships that may encode root configuration patterns.\\


Using the decomposition method (8 features), we count the number of near-zero coefficients (absolute value $< 0.1$), the ratio of maximum to minimum coefficient magnitudes, coefficient variance, an indicator for potential $x$-factorization (whether $|E| < 0.1$), and four normalized measures testing for factorization patterns inspired by Sylvester matrix relationships.\\

\subsection{Distillation Protocol}

To interpret the nonlinear decision boundary learned by neural networks, we employ a knowledge distillation framework to extract an interpretable surrogate decision tree. We first train the neural network with the identified features. We then use the trained network to generate predictions on the training data and train the decision tree to mimic them. After the tree is trained, we evaluate it on a test set, calculating the tree's (1) fidelity (agreement with network predictions on unseen data), and (2) balanced accuracy (agreement with ground truth). Feature importance is also analyzed via SHAP, which is SHapley Additive exPlanations \cite{lundberg2017}.

\section{Results and Analysis}

\subsection{Validation of Methodology}

For quadratic polynomials, neural networks achieved 99.3\% $\pm$ 0.1\% accuracy using raw coefficients, rising to perfect 100.0\% $\pm$ 0.0\% when provided the discriminant ratio $b^2/ac$. Decision trees achieved perfect 100\% classification using this single feature, learning threshold values near 4.0 that match theoretical predictions. For cubic polynomials, neural networks achieved 98.5\% $\pm$ 0.002\% accuracy using raw coefficients, while decision trees, using the invariant ratio $\beta^2/\alpha^3$ returned the threshold 3.999, achieving perfect 100\% accuracy using a threshold nearly identical to the theoretical value of 4.0 separating real from complex root cases. For quartic polynomials, neural networks achieved 97.5\% $\pm$ 0.2\% balanced accuracy with raw coefficients. When provided classical invariants (discriminant expressions and Tschirnhaus-style features), decision trees reached 99.4\% and other ensemble methods (specifically Random Forest and Gradient Boosting) achieved 100.0\% balanced accuracy. These validation results establish that our methodology can extract meaningful mathematical structure when ground truth is available and appropriate features are provided, giving us confidence to apply the same approach to quintics for which no closed-form solutions are known. All symbolic regression experiments were conducted with fixed search depth, operator set, and computational budget. While increasing these resources may improve performance, our results indicate that under practical search constraints, symbolic regression struggles beyond quadratic structure. The Symbolic Regression results reveal a distinct ``Complexity Cliff,'' demonstrating how performance degrades as polynomial complexity increases, as summarized in Table \ref{tab:complexity_cliff}.

\begin{table}[htbp]
\centering
\caption{Summary of Symbolic Regression Performance by Polynomial Degree}
\label{tab:complexity_cliff}
\begin{tabular}{c c l}
\hline
\textbf{Polynomial} & \textbf{Accuracy (\%)} & \textbf{Observations} \\
\hline
Quadratic & $97.3 \pm 4.2$ & Discovered discriminant rule (100\% in 2/5 trials) \\
Cubic & $94.4 \pm 1.6$ & Discovered highly accurate but complex heuristic \\
Quartic & $47.3 \pm 5.8$ & Consistently failed; relied on simple linear heuristic \\
Quintic & $53.5 \pm 2.9$ & Failed to find interpretable rules; relied on complex heuristics \\
\hline
\end{tabular}
\end{table}

This confirms that while Symbolic Regression can discover simple algebraic rules, classification problems beyond the second degree pose a complexity barrier that prevents autonomous symbolic discovery, in contrast to the 84\%+ performance of neural networks.

\subsection{Model Screening Results}

In our initial screening phase, we evaluated all model classes on quintic polynomial classification using raw coefficients across three separate seeds. As shown in Table \ref{tab:screening}, Neural networks consistently achieved the highest performance, while interpretable models struggled significantly to establish boundaries without guided features.

\begin{table}[htbp]
\caption{Degree 5 Initial Screening Performance (Mean $\pm$ 95\% CI)}
\centering
\renewcommand{\arraystretch}{1.2}
\begin{tabular}{lc}
\toprule
\textbf{Model} & \textbf{Balanced Accuracy (Raw Coefficients)} \\
\midrule
Neural Networks & 83.8\% $\pm$ 3.9\% \\
Gradient Boosting & 63.3\% $\pm$ 3.8\% \\
SVM (RBF) & 62.9\% $\pm$ 0.4\% \\
Random Forest & 62.5\% $\pm$ 1.5\% \\
Decision Trees (CART) & 53.1\% $\pm$ 3.1\% \\
Logistic Regression & 41.6\% $\pm$ 0.9\% \\
\bottomrule
\end{tabular}
\label{tab:screening}
\end{table}
 
Based on these findings, we focused our rigorous 20-trial statistical analysis on neural networks and decision trees, as these two model classes directly address our research question about learning patterns and extracting interpretable rules.
As shown in Table \ref{t4}, testing neural networks and decision trees across 20 independent trials with different random seeds reveals a significant gap between their capabilities when using only raw polynomial coefficients. The inclusion of the \textit{Crit8} feature (sign changes at critical points) significantly narrows the capability gap between the black-box neural model and the interpretable decision tree (84.2\% versus 89.9\%).

\begin{table}[htbp]
    \centering
    \begin{threeparttable}
        \caption{Performance comparison of neural networks and decision trees when additional feature ``crit8'' added.}
        \label{t4}
        \begin{tabular}{lcc}
            \toprule
            \textbf{Model} & \textbf{Raw Coefficients} & \textbf{With Crit8\tnote{1}} \\
            \midrule
            Neural Networks & 84.3\% $\pm$ 0.9\% & 89.9\% $\pm$ 1.4\% \\
            Decision Trees & 59.9\% $\pm$ 0.9\% & 84.2\% $\pm$ 1.2\% \\
            \midrule
            \textbf{Performance Gap} & \textbf{24.4 points} & \textbf{5.7 points} \\
            \bottomrule
        \end{tabular}
        \begin{tablenotes}
            \small
            \item[1] Crit8 Feature: the number of sign changes in the sequence $\{p(c_1), p(c_2), \dots, p(c_k)\}$ where $c_i$ are ordered critical points.
        \end{tablenotes}
    \end{threeparttable}
\end{table}

The 33\% baseline reflects random guessing among three classes. Neural networks achieved 84.3\% balanced accuracy using only raw coefficients, substantially above baseline. However, decision trees achieved only 59.9\%, demonstrating their inability to autonomously discover useful classification rules from raw data. While the performance of Neural Networks was impressive, the performance gap highlights a fundamental limitation. The structure that neural networks learn internally cannot be extracted as interpretable rules without additional guidance. When provided Crit8 (sign changes at critical points) as an explicit feature, decision trees improve dramatically to 84.2\% (+24.3 points). Neural networks improve modestly (+5.6 points). This asymmetry suggests that the neural network potentially learned a rule similar to Crit8 with raw coefficients alone, while interpretable models require human-engineered features to achieve good performance on this problem.\\

In an attempt to better understand the decision boundary learned by neural networks, we trained them on all 63 mathematical features (see subsection 2.2) for quintic classification, then attempted to distill their knowledge into decision trees. As shown in Table \ref{tab:distillation_results}, the 98.9\% test-set fidelity demonstrates that decision trees can successfully mimic neural network predictions when provided with engineered features. The 84.6\% standalone balanced accuracy confirms that the extracted rules are useful.\\

\begin{table}[htbp]
    \centering
    \begin{threeparttable}
        \caption{Quantitative results for the Knowledge Distillation process from Neural Networks (NN) to Decision Trees. The high fidelity indicates that the interpretable tree successfully approximates the neural model's decision logic.}
        \label{tab:distillation_results}
        \begin{tabular}{lll}
            \toprule
            \textbf{Metric} & \textbf{Value} & \textbf{Interpretation} \\
            \midrule
            NN Test Balanced Accuracy & 87.0\% $\pm$ 1.2\% & Robust feature learning \\
            \textbf{Tree Test Fidelity} & \textbf{98.9\% $\pm$ 0.1\%} & Successful logic transfer \\
            Tree Standalone Accuracy & 84.6\% $\pm$ 1.2\% & Generalizable rule extraction \\
            \bottomrule
        \end{tabular}
    \end{threeparttable}
\end{table}

The SHAP analysis \cite{lundberg2017}  revealed that the feature Crit8, which counted sign changes in critical points, enabled this success. It boasted a feature importance score of 97.5\%, 145.6$\times$ more important than the second-ranked feature. Examining the best-trial decision tree revealed the classification structure as follows.

\begin{equation*}
\text{Classification} = \begin{cases}
\text{1 real, 4 complex} & \text{if Crit8} \leq 0.5 \\
\text{3 real, 2 complex} & \text{if } 0.5 < \text{Crit8} \leq 1.5 \\
\text{5 real roots} & \text{if Crit8} > 1.5
\end{cases}
\end{equation*}

These rules have a clear mathematical foundation. Critical points where $p'(x) = 0$ divide the real line into monotonic regions. By the Intermediate Value Theorem, sign changes in $\{p(c_1), \ldots, p(c_k)\}$ indicate zero crossings. Polynomials with 5 real roots exhibit multiple oscillations (high Crit8 $\geq$ 2); those with 1 real root show near-monotonic behavior (low Crit8 $\leq$ 1). The learned thresholds (0.5, 1.5) match this mathematical intuition. However, it is crucial to note that this interpretable rule was only discoverable because we explicitly engineered the Crit8 feature. The neural network did not produce this rule autonomously, and without our feature engineering, the decision tree could not have found it. Additionally, testing each classical method individually yields these results (see Table \ref{tab:feature_ablation}). Only when critical point features are explicitly included does neural network performance improve above baseline; other methods often degrade it through noise or overfitting. For decision trees, critical points provide massive improvement (+24.3 points) while other methods offer minimal gains.

\begin{table}[htbp]
    \centering
    \begin{threeparttable}
        \caption{Comparative performance across different mathematical feature sets for quintic polynomials. Bold values indicate the optimal feature configuration, where \textit{Critical points} provide the most significant gain in interpretability and accuracy.}
        \label{tab:feature_ablation}
        \begin{tabular}{lcc}
            \toprule
            \textbf{Feature Set} & \textbf{Neural Networks} & \textbf{Decision Trees} \\
            \midrule
            Raw coefficients & 84.3\% $\pm$ 0.9\% & 59.9\% $\pm$ 0.9\% \\
            + Sturm sequences & 85.2\% $\pm$ 1.1\% & 62.1\% $\pm$ 4.1\% \\
            + Newton sums & 83.0\% $\pm$ 1.3\% & 63.0\% $\pm$ 4.2\% \\
            + \textbf{Critical points} & \textbf{89.3\% $\pm$ 0.7\%} & \textbf{84.2\% $\pm$ 1.2\%} \\
            + Hybrid symbolic & 80.5\% $\pm$ 1.4\% & 61.4\% $\pm$ 4.0\% \\
            \midrule
            All 63 combined & 87.0\% $\pm$ 1.2\% & 68.2\% $\pm$ 4.5\% \\
            \bottomrule
        \end{tabular}
    \end{threeparttable}
\end{table}

\subsection{Model Robustness and Generalization}

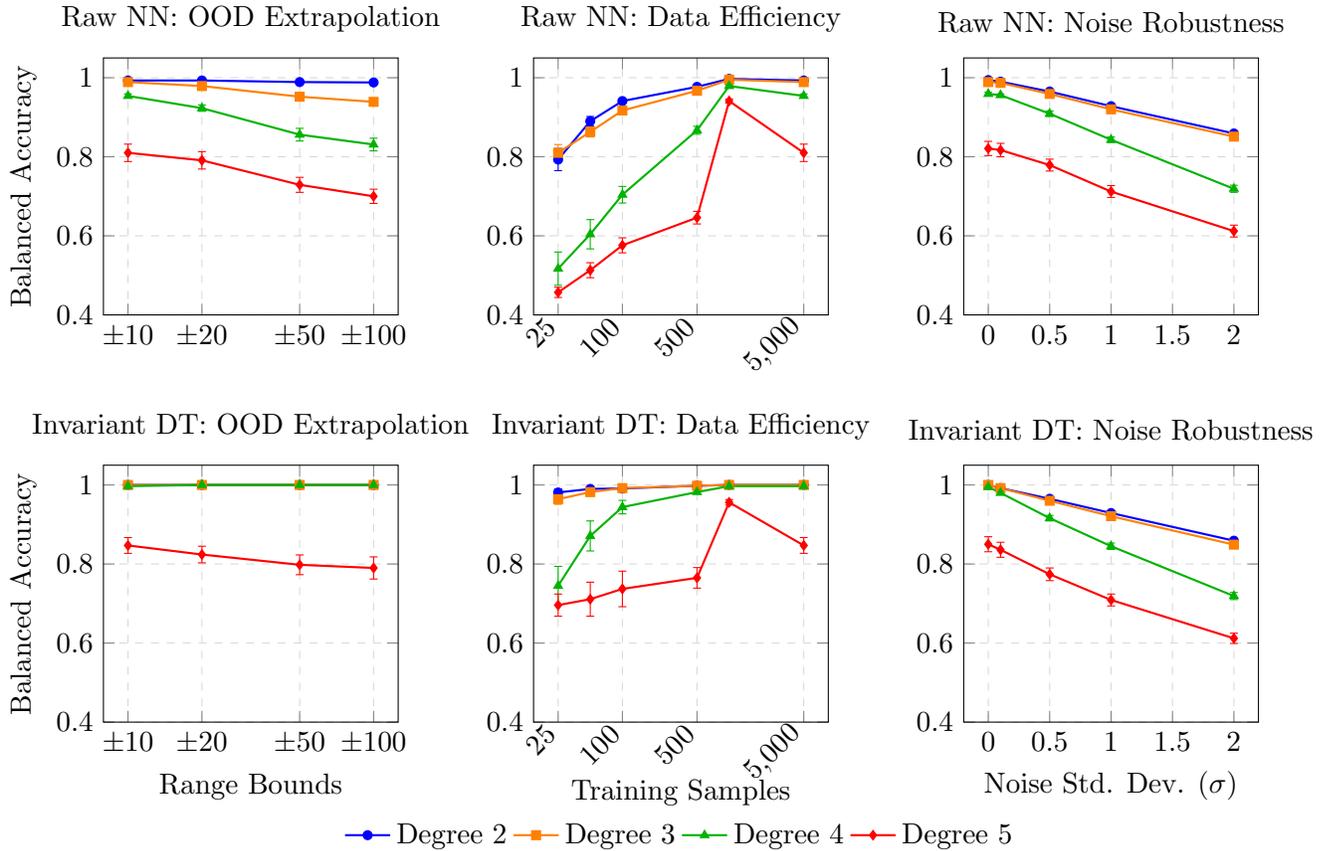
\begin{figure}[htbp]
    \centering
    \begin{tikzpicture}
        \begin{groupplot}[
            group style={
                group size=3 by 2,
                horizontal sep=1.8cm,
                vertical sep=2.0cm,
            },
            width=5.5cm,
            height=5.0cm,
            ymin=0.4, ymax=1.05,
            grid=major,
            grid style={dashed, gray!30},
            ylabel style={align=center},
            xlabel style={align=center},
            every axis plot/.append style={thick, mark size=1.5pt},
            error bars/y dir=both,
            error bars/y explicit
        ]
        
        
        \nextgroupplot[
            xmode=log,
            title={Raw NN: OOD Extrapolation},
            ylabel={Balanced Accuracy},
            xtick={10, 20, 50, 100},
            xticklabels={$\pm 10$, $\pm 20$, $\pm 50$, $\pm 100$},
            log ticks with fixed point
        ]
        \addplot[color=blue, mark=*] coordinates {(10,0.993) +- (0,0.001) (20,0.993) +- (0,0.002) (50,0.989) +- (0,0.003) (100,0.988) +- (0,0.003)};
        \addplot[color=orange, mark=square*] coordinates {(10,0.989) +- (0,0.001) (20,0.979) +- (0,0.002) (50,0.952) +- (0,0.004) (100,0.939) +- (0,0.005)};
        \addplot[color=green!70!black, mark=triangle*] coordinates {(10,0.954) +- (0,0.005) (20,0.923) +- (0,0.008) (50,0.856) +- (0,0.016) (100,0.831) +- (0,0.016)};
        \addplot[color=red, mark=diamond*] coordinates {(10,0.810) +- (0,0.022) (20,0.791) +- (0,0.022) (50,0.729) +- (0,0.019) (100,0.700) +- (0,0.018)};
        
        \nextgroupplot[
            xmode=log,
            title={Raw NN: Data Efficiency},
            xtick={25, 100, 500, 5000},
            log ticks with fixed point,
            x tick label style={rotate=45, anchor=east}
        ]
        \addplot[color=blue, mark=*] coordinates {(25,0.793) +- (0,0.028) (50,0.890) +- (0,0.012) (100,0.941) +- (0,0.005) (500,0.977) +- (0,0.002) (1000,0.997) +- (0,0.001) (5000,0.993) +- (0,0.001)};
        \addplot[color=orange, mark=square*] coordinates {(25,0.810) +- (0,0.021) (50,0.863) +- (0,0.013) (100,0.917) +- (0,0.009) (500,0.967) +- (0,0.003) (1000,0.995) +- (0,0.001) (5000,0.989) +- (0,0.001)};
        \addplot[color=green!70!black, mark=triangle*] coordinates {(25,0.517) +- (0,0.042) (50,0.604) +- (0,0.037) (100,0.704) +- (0,0.021) (500,0.867) +- (0,0.010) (1000,0.979) +- (0,0.002) (5000,0.954) +- (0,0.005)};
        \addplot[color=red, mark=diamond*] coordinates {(25,0.457) +- (0,0.013) (50,0.513) +- (0,0.019) (100,0.576) +- (0,0.019) (500,0.646) +- (0,0.016) (1000,0.941) +- (0,0.005) (5000,0.810) +- (0,0.022)};
        
        \nextgroupplot[
            title={Raw NN: Noise Robustness},
            xtick={0, 0.5, 1.0, 1.5, 2.0}
        ]
        \addplot[color=blue, mark=*] coordinates {(0,0.994) +- (0,0.001) (0.1,0.990) +- (0,0.002) (0.5,0.965) +- (0,0.003) (1.0,0.928) +- (0,0.004) (2.0,0.859) +- (0,0.005)};
        \addplot[color=orange, mark=square*] coordinates {(0,0.989) +- (0,0.001) (0.1,0.987) +- (0,0.002) (0.5,0.959) +- (0,0.004) (1.0,0.920) +- (0,0.003) (2.0,0.851) +- (0,0.005)};
        \addplot[color=green!70!black, mark=triangle*] coordinates {(0,0.959) +- (0,0.004) (0.1,0.956) +- (0,0.003) (0.5,0.909) +- (0,0.007) (1.0,0.843) +- (0,0.007) (2.0,0.719) +- (0,0.009)};
        \addplot[color=red, mark=diamond*] coordinates {(0,0.821) +- (0,0.018) (0.1,0.817) +- (0,0.017) (0.5,0.779) +- (0,0.015) (1.0,0.712) +- (0,0.015) (2.0,0.612) +- (0,0.015)};


        \nextgroupplot[
            xmode=log,
            title={Invariant DT: OOD Extrapolation},
            xlabel={Range Bounds},
            ylabel={Balanced Accuracy},
            xtick={10, 20, 50, 100},
            xticklabels={$\pm 10$, $\pm 20$, $\pm 50$, $\pm 100$},
            log ticks with fixed point
        ]
        \addplot[color=blue, mark=*] coordinates {(10,1.000) +- (0,0) (20,1.000) +- (0,0) (50,1.000) +- (0,0) (100,1.000) +- (0,0)};
        \addplot[color=orange, mark=square*] coordinates {(10,1.000) +- (0,0) (20,1.000) +- (0,0) (50,1.000) +- (0,0) (100,1.000) +- (0,0)};
        \addplot[color=green!70!black, mark=triangle*] coordinates {(10,0.997) +- (0,0.001) (20,1.000) +- (0,0) (50,1.000) +- (0,0) (100,1.000) +- (0,0)};
        \addplot[color=red, mark=diamond*] coordinates {(10,0.847) +- (0,0.020) (20,0.824) +- (0,0.021) (50,0.798) +- (0,0.025) (100,0.790) +- (0,0.028)};

        \nextgroupplot[
            xmode=log,
            title={Invariant DT: Data Efficiency},
            xlabel={Training Samples},
            xlabel style={yshift=11pt}, 
            xtick={25, 100, 500, 5000},
            log ticks with fixed point,
            x tick label style={rotate=45, anchor=east},
            legend style={at={(0.5,-0.35)}, anchor=north, legend columns=4, draw=none, fill=none}
        ]
        \addplot[color=blue, mark=*] coordinates {(25,0.981) +- (0,0.009) (50,0.990) +- (0,0.003) (100,0.991) +- (0,0.005) (500,0.998) +- (0,0.001) (1000,1.000) +- (0,0) (5000,1.000) +- (0,0)};
        \addplot[color=orange, mark=square*] coordinates {(25,0.964) +- (0,0.013) (50,0.982) +- (0,0.010) (100,0.992) +- (0,0.003) (500,0.998) +- (0,0.001) (1000,1.000) +- (0,0) (5000,1.000) +- (0,0)};
        \addplot[color=green!70!black, mark=triangle*] coordinates {(25,0.745) +- (0,0.049) (50,0.871) +- (0,0.038) (100,0.944) +- (0,0.017) (500,0.982) +- (0,0.004) (1000,0.997) +- (0,0.001) (5000,0.997) +- (0,0.001)};
        \addplot[color=red, mark=diamond*] coordinates {(25,0.696) +- (0,0.028) (50,0.711) +- (0,0.043) (100,0.737) +- (0,0.045) (500,0.765) +- (0,0.026) (1000,0.956) +- (0,0.007) (5000,0.847) +- (0,0.020)};
        \vspace{1.2cm}
        \addlegendentry{Degree 2}
        \addlegendentry{Degree 3}
        \addlegendentry{Degree 4}
        \addlegendentry{Degree 5}

        \nextgroupplot[
            title={Invariant DT: Noise Robustness},
            xlabel={Noise Std. Dev. ($\sigma$)},
            xtick={0, 0.5, 1.0, 1.5, 2.0}
        ]
        \addplot[color=blue, mark=*] coordinates {(0,1.000) +- (0,0) (0.1,0.992) +- (0,0.002) (0.5,0.965) +- (0,0.003) (1.0,0.929) +- (0,0.004) (2.0,0.859) +- (0,0.005)};
        \addplot[color=orange, mark=square*] coordinates {(0,1.000) +- (0,0) (0.1,0.992) +- (0,0.001) (0.5,0.960) +- (0,0.004) (1.0,0.921) +- (0,0.003) (2.0,0.849) +- (0,0.005)};
        \addplot[color=green!70!black, mark=triangle*] coordinates {(0,0.995) +- (0,0.002) (0.1,0.980) +- (0,0.004) (0.5,0.916) +- (0,0.007) (1.0,0.845) +- (0,0.008) (2.0,0.719) +- (0,0.009)};
        \addplot[color=red, mark=diamond*] coordinates {(0,0.850) +- (0,0.019) (0.1,0.836) +- (0,0.019) (0.5,0.774) +- (0,0.016) (1.0,0.709) +- (0,0.015) (2.0,0.612) +- (0,0.013)};

        \end{groupplot}
    \end{tikzpicture}
    \caption{Behavioral divergence between neural approximation (Top Row) and discrete mathematical execution (Bottom Row). Error bars indicate 95\% Confidence Intervals. Left: Invariant decision trees maintain perfectly stable, scale-invariant bounds out-of-distribution, whereas neural approximations decay. Middle: Trees lock in boundaries with mere dozens of samples, while networks require thousands. Right: Both models degrade identically under noise, reflecting the fundamental alteration of the true polynomial roots.}
    \label{fig:robustness}
\end{figure}

To rigorously evaluate the nature of the neural network's internal representations, we established a mathematically ``desired'' control group. We compared the performance of black-box neural networks trained on raw coefficients against interpretable decision trees equipped with exact algebraic invariants (e.g., the exact discriminants $\Delta$ for lower degrees, and the Crit8 features for degree 5). Both models were subjected to out-of-distribution (OOD) extrapolation, data efficiency, and input noise robustness stress tests using 5-Fold Stratified Cross-Validation across 3 random seeds. The comparative performance is visualized in Figure \ref{fig:robustness}.

The stress tests provide strong empirical evidence that the neural network did not recover discrete algebraic invariants in a scale-invariant or symbolically extractable form. A true mathematical invariant (such as $\Delta = b^2 - 4ac$) is perfectly scale-invariant and requires virtually no training data to define a decision boundary. When evaluated out-of-distribution on coefficients scaled up to $10\times$ their training bounds, the control group (invariant decision trees) maintained flawless stability, scoring an exact $1.000 \pm 0.000$ balanced accuracy across degrees 2, 3, and 4. In stark contrast, the raw neural networks failed to generalize their approximations; for example, quartic accuracy dropped from $95.4\% \pm 0.5\%$ to $83.1\% \pm 1.6\%$. Furthermore, the data efficiency curves reveal the underlying learning mechanisms. The invariant control models achieved near-perfect accuracy with as few as 50 to 100 samples. Conversely, the neural networks required thousands of samples to navigate the complex decision geometries of higher-degree polynomials, with quintic models continuing to improve steadily up to the 5,000 sample maximum without plateauing perfectly. Interestingly, both the raw neural networks and the invariant control models degraded identically under continuous Gaussian noise (e.g., degree 4 accuracy dropping to roughly $71.9\%$ at $\sigma=2.0$ for both models). This synchronous failure highlights the chaotic sensitivity of the polynomial coefficient space itself: introducing noise fundamentally alters the roots of the polynomial, meaning even perfect mathematical logic will inherently misclassify the original data labels. Ultimately, these side-by-side results clearly illustrate the interpretability gap. While the neural network finds a highly accurate, resilient representation of the math, it does so by building a continuous, data-dependent, localized geometric approximation bounded by its training space. Because the knowledge is stored geometrically rather than symbolically, it cannot be easily extracted. Translating that latent pattern-matching capability into explicit, human-readable logic strictly requires human intervention and guided feature engineering. These results suggest that the network’s learned representation is geometric and data-dependent rather than symbolic and invariant. Because the knowledge is stored geometrically rather than symbolically, it cannot be easily extracted as explicit mathematical logic without human intervention.

\section{Discussion and Conclusion}

Our extensive evaluation of ML methods for quintic polynomial classification yields a nuanced conclusion. While neural networks can achieve reasonable balanced accuracy (84.3\%) from raw coefficients, the models we tested did not autonomously discover interpretable mathematical rules as we had hoped. Interpretable models like decision trees perform poorly on raw data (59.9\%), and extracting rules from neural networks requires explicit human-guided feature engineering. When we provided the Crit8 feature, decision trees achieved 84.2\% balanced accuracy and produced interpretable rules. But this required us to identify and implement the relevant mathematical concept; the models did not discover it on their own.\\

The present study carries a number of limitations. First, our experiments rely on polynomials with coefficients sampled uniformly from $[-10, 10],$ which defines a specific coefficient geometry; alternative sampling distributions may alter root configuration frequencies and learned decision boundaries. Second, root classification is determined numerically using a fixed tolerance threshold, meaning near-multiple roots may introduce unavoidable labeling instability. Third, while we evaluated a broad range of classical and modern ML models, including neural networks and symbolic regression, our analysis does not exhaust the space of possible architectures (e.g., transformer-based symbolic learners, equivariant networks, or program-synthesis systems). Fourth, symbolic regression experiments were conducted under practical computational constraints, and deeper search regimes may yield different outcomes. Finally, although our stress tests strongly suggest that neural networks learn continuous geometric approximations rather than discrete symbolic invariants, we did not directly probe internal representations; therefore, we cannot exclude the possibility that latent algebraic structure exists in a form not accessible through distillation. These considerations suggest that our conclusions should be interpreted within the scope of the evaluated ML models and experimental design.\\

Our results do not prove that no ML approach can autonomously discover interpretable mathematical rules. However, they suggest this remains an open challenge. Future work may explore richer symbolic regression systems, transformer architectures trained on formal mathematics, hybrid neuro-symbolic methods, and alternative interpretability strategies beyond distillation. In conclusion, bridging the gap between geometric approximation and symbolic reasoning remains a central challenge for interpretable machine learning in structured mathematical domains.

\section*{Declaration of competing interest}

\noindent The authors declare no conflicts of interest.

\section*{Data Availability}

\noindent All codes to perform numerical simulations, model fitting and hypothesis testing are available in our \href{https://github.com/rohanthomas334-sketch/MLPolyRoots}{GitHub repository}. The README file contains instructions to use each code to verify the numerical results.

\bibliographystyle{unsrtnat}
\bibliography{references}

\end{document}